\documentclass[12pt]{article}
\usepackage{amssymb}
\usepackage{amsmath, amscd, amsthm}
\usepackage{amsrefs}
\usepackage{mathrsfs}
\usepackage[all]{xy}
\usepackage{amsfonts}

\begin{document}

\theoremstyle{remark}
\newtheorem*{Remark}{Remark}
\newtheorem*{Remarks}{Remarks}
\newtheorem*{Proof}{Proof}
\theoremstyle{definition}
\newtheorem{Definition}{Definition}
\newtheorem*{Example}{Example}
\newtheorem*{Examples}{Examples}

\theoremstyle{plain}
\newtheorem{Theorem}{Theorem}
\newtheorem{Theorema}{Theorem}
\newtheorem{Proposition}{Proposition}
\newtheorem{Corollary}{Corollary}
\newtheorem{Lemma}{Lemma}

\newcommand{\g}{\mathfrak g}
\newcommand{\cl}{\mathbf{Cl}}
\newcommand{\vs}{\vspace{5mm}}
\renewcommand{\labelenumi}{\roman{enumi})}
\setlength{\parindent}{0em}
\small\normalsize

\begin{center}
{\Large\sc On a twisted Reidemeister torsion 
\footnote{\today}\footnote{Mathematics Subject Classification (2010): 57Q10}
}\\ $\ $ \\
{ Ricardo Garc\'{\i}a L\'opez
\footnote{
 Departament d'\`Algebra i Geometria.
 Universitat de Barcelona,
 Gran Via, 585.
 E-08007, Barcelona, Spain.
 e-mail: ricardogarcia@ub.edu}
}
\end{center}
\begin{abstract}
Given a finite simplicial complex, a unimodular representation of its fundamental group and a closed twisted cochain of odd degree, we define a twisted version of the Reidemeister torsion, extending a previous definition of V. Mathai and S. Wu. 
The main tool is a complex of piecewise smooth currents introduced by J. Dupont in \cite{Du}.
\end{abstract}

{\bf  Introduction.}\vs

Let $X$ be a compact, oriented, smooth manifold, denote $\Omega^{\bullet}(X)$ the de Rham complex of smooth forms on $X$. Let $(\mathcal E, \nabla)$ be a vector bundle with flat connection over $X$, denote $\Omega^{\bullet}(X, \mathcal E) = (\Omega^{\bullet}\otimes \mathcal E)(X)$ the $\mathcal E$--valued differential forms on $X$, recall that the connection $\nabla: \mathcal E(X)\longrightarrow \Omega^1(X,\mathcal E)$ extends to a  differential $\delta: \Omega^{\bullet}(X, \mathcal E)\longrightarrow \Omega^{\bullet+1}(X, \mathcal E)$.
If
$T\in \Omega^{\bullet}(X)$ denotes a closed, odd--degree differential form, then the twisted de Rham complex attached to these data is the $\mathbb Z_2$--graded\footnote{Here, $\mathbb Z_2:=\mathbb Z/ 2\mathbb Z$} complex $((\Omega^{\bullet}\otimes \mathcal E)(X), \delta + T\,\wedge\,)$. When $T$ is of degree three, its cohomology was first considered by R. Rohm and E. Witten in \cite{RW}, in connection with string theory. It is the recipient of the Chern character in twisted $K$-theory, as defined by P. Bouwknegt et. al. in \cite{BCM} or by M. Atiyah and G. Segal in \cite{AS}.
\vs

A simplicial counterpart of the twisted de Rham cohomology was considered by V. Mathai and S. Wu in \cite{MW}: Let $K$ be a finite simplicial complex,
$\rho: \pi_1(K) \longrightarrow GL(E)$ a  representation of its fundamental group, let $(C^{\bullet}(K,E), \partial)$ denote the complex of simplicial cochains twisted by $\rho$. If $\vartheta\in \oplus_{i\geqslant 1}C^{2i+1}(K)$ is a cocycle
of degree greater than $\dim K/2$, 
they consider the  $\mathbb Z_2$--graded complex $(C^{\bullet}(K,E), \partial +\vartheta\,\cup)$, its cohomology $H^{even}(K,E, \vartheta)$ and $H^{odd}(K,E, \vartheta)$ is the  twisted simplicial cohomology attached to $K,E,\vartheta$\footnote{In a sense, this is a doubly twisted cohomology, by the representation $\rho$ and by the cocycle $\vartheta$.}. When the representation $\rho$ is unimodular, Mathai and Wu consider also a twisted version of the Reidemeister torsion, which is an element of the determinant line $\det H^{\bullet}(K,E, \vartheta)$. If $K$ is a smooth triangulation of a smooth manifold and $\vartheta$ is of 
maximal degree, then it maps under a canonical isomorphism to the twisted version of the Ray-Singer torsion. The last 
lives in the 
determinant line of the twisted de Rham cohomology, and it is the main object of study in \cite{MW}. \vs

In this note, using a complex of simplicial currents introduced by J.L. Dupont in \cite{Du}, we show that one can  define a twisted simplicial cohomology and a twisted Reidemeister torsion 
without the restriction $\deg \vartheta > \dim K/2$. For smooth manifolds, this twisted simplicial cohomology  is  isomorphic to the twisted de Rham cohomology. Also, the torsion defined here maps under a canonical isomorphism of determinant lines to the twisted Reidemeister torsion introduced by Mathai and Wu, when the later is defined. Invariance under subdivision holds both for the twisted simplicial cohomology and for the twisted torsion. If the representation $\rho$ is trivial then, instead of Dupont's currents one can define twisted cohomology just using the Sullivan-Thom-Whitney algebra of piecewise smooth differential forms, this was indeed done by D. Freed et al. in \cite{FHTe}*{(2.1)}. \vs

I thank the anonymous referee for her/his comments. \vs

{\bf  Twisted cohomology and twisted Reidemeister torsion.}\vs

{\it Determinants\ }{(see \cite{KM}, \cite{K}):} Let $\mathbf k$ be a field of characteristic zero, let $Pic^{\, \mathbb Z}(\mathbf k)$ denote the category of graded $\mathbf k$--lines. An object of
$Pic^{\mathbb Z}(\mathbf k)$ is a pair $L=(\ell, n)$ where $\ell$ is a $1$-dimensional $\mathbf k$-vector space and
$n\in\mathbb Z$. In this category, $Hom((\ell_1,n_1), (\ell_2,n_2))=\emptyset$ if $n_1\neq n_2$
and $Hom((\ell_1,n), (\ell_2,n))=Hom_{\mathbf k}(\ell_1, \ell_2)-\{0\}$. 
There is a tensor product functor $\otimes: Pic^{\,\mathbb Z}(\mathbf k)\times  Pic^{\,\mathbb Z}(\mathbf k) \longrightarrow Pic^{\,\mathbb Z}(\mathbf k)$, defined on objects by
\[
((\ell_1,n_1), (\ell_2,n_2)) \longmapsto (\ell_1\otimes \ell_2, n_1+n_2)\,.
\]
It is proved in loci cit. that $Pic^{\mathbb Z}(\mathbf k)$ is a Picard category, with unit $(\mathbf k, 0)$, the standard associativity constraint and the commutativity constraint 
\[
(\ell_1,n_1)\otimes (\ell_2,n_2) \stackrel{\sim}\longrightarrow (\ell_2,n_2)\otimes (\ell_1,n_1)
\]
given by $u\otimes w \longmapsto (-1)^{n_1n_2} w\otimes u$. In particular, given a family of graded lines $\{L_i\}_{i\in I}$, the tensor product 
$\otimes_{i\in I}L_i$ is unambiguously defined up to a canonical, functorial isomorphism. If $L_1, L_2$ are graded lines, we write $L_1=L_2$ if they are canonically isomorphic. 
If $L=(\ell, n)$ is a graded line, we denote $L^{-1}=(Hom_{\mathbf k}(\ell,\mathbf k), -n)$.
For $(\ell, n)$ a graded line, the reference to the integer $n$ will be omitted when no confusion arises. \vs

Given a $\mathbf k$--vector space $V$ of dimension $n>0$, put $\det V= (\bigwedge^n V, n)$ (and put $\det\{0\}=(\mathbf k, 0)$). Let
\[
\mathscr{C}: \ \ \  0 \longrightarrow \mathscr C^r \stackrel{\partial_r}{\longrightarrow} \mathscr C^{r+1} \stackrel{\partial_{r+1}}{\longrightarrow}... \longrightarrow \mathscr C^s \longrightarrow 0 
\]
be a bounded complex of finite dimensional $\mathbf k$--vector spaces. The determinant of $\mathscr{C}$ is the graded line 
\[
\det \mathscr{C}:= \bigotimes_{i=r}^s (\det \mathscr C^i)^{(-1)^{i}}.
\]
Regarding the cohomology of $\mathscr{C}$ as a complex with zero differential, we can consider also 
\[
 \det H^{\bullet}(\mathscr{C}, \partial)=\bigotimes_{i=r}^s (\det H^{i}(\mathscr{C}, \partial))^{(-1)^{i}}\,,
\]
and it is proved in \cite{KM} that there is a canonical, functorial isomorphism
$
\det \mathscr{C}\longrightarrow \det H^{\bullet}(\mathscr{C}, \partial)
$
that we call the Knudsen-Mumford isomorphism.\vs

{\it Filtered $\mathbb Z_2$-graded complexes:} Let $\mathbf C$ be $\mathbb Z_2$-graded 
$\mathbf k$-vector space, that is, a $\mathbf k$-vector space endowed with a decomposition 
$\mathbf C=\mathbf C^{ev}\oplus \mathbf C^{od}$.  
A $\mathbb Z_2$-graded complex is a couple $(\mathbf C,\partial)$ where 
$\partial: \mathbf C \longrightarrow \mathbf C$ is a $\mathbf k$-linear map such that $\partial^2=0$, $\partial(\mathbf C^{ev})\subset \mathbf C^{od}$, $\partial(\mathbf C^{od})\subset \mathbf C^{ev}$. A morphism of $\mathbb Z_2$-graded complexes is a linear map preserving both the differentials and the decompositions. The twisted cohomology of $(\mathbf C,\partial)$ is defined by
\begin{eqnarray*}
H^{ev}(\mathbf C, \partial) &=& \frac{\mbox{ker}(\partial: \mathbf C^{ev}\longrightarrow \mathbf C^{od} )}{\mbox{im} (\partial: \mathbf C^{od}\longrightarrow \mathbf C^{ev})} \,,\hspace{0.2cm}\\
H^{od}(\mathbf C, \partial)  &=& \frac{\mbox{ker}(\partial: \mathbf C^{od}\longrightarrow \mathbf C^{ev}) }{\mbox{im} (\partial: \mathbf C^{ev}\longrightarrow \mathbf C^{od})}\,.
\end{eqnarray*}
If $\mathbf C$ is finite dimensional, put
\begin{eqnarray*}
 \det \mathbf C &=& \det \mathbf C^{ev} \otimes_{\mathbf k} (\det \mathbf C^{od})^{-1} \,,\hspace{0.2cm}\\
 \det H^{\bullet} (\mathbf C, \partial) &=& \det H^{ev}(\mathbf C, \partial) \otimes_{\mathbf k} (\det H^{od}(\mathbf C, \partial))^{-1}\,.
\end{eqnarray*}\vspace{0.5mm}

\begin{Lemma} 
If $(\mathbf C, \partial)$ is a $\mathbb Z_2$-graded complex and $\mathbf C$ is finite dimensional, then $ \det \mathbf C= \det H^{\bullet} (\mathbf C, \partial)$.
\end{Lemma}
\begin{Proof}
Put $A=\mbox{im}(\partial: \mathbf C^{od}\longrightarrow \mathbf C^{ev})$, $B=\mbox{ker} (\partial: \mathbf C^{od}\longrightarrow \mathbf C^{ev})$. Consider the complex
\[
\stackrel{\bullet}0\,\longrightarrow A \stackrel{i}\longrightarrow \mathbf C^{ev}\stackrel{\partial}\longrightarrow \mathbf C^{od} \stackrel{\pi}\longrightarrow \mathbf C^{od}/B \longrightarrow 0
\]
where the dot denotes the term placed at zero, $i$ is the inclusion and $\pi$ the projection. The Knudsen-Mumford isomorphism  gives
\[
\det H^{ev}(\mathbf C, \partial) \otimes \det H^{od}(\mathbf C, \partial)^{-1} = \det A ^{-1}\otimes \det \mathbf{C}^{ev}\otimes (\det \mathbf{C}^{od})^{-1}
\otimes \det \mathbf{C}^{od}/B
\]
and since $\mathbf C^{od}/B \cong A$, the right hand side is canonically isomorphic to $\det \mathbf C$, so we are done. $\Box$
\end{Proof}\vs

If ($\mathbf C, \partial)$ is a $\mathbb Z_2$-graded complex, we denote $\widetilde{\mathbf C}$ the ($\mathbb Z$-graded) complex defined by
\begin{equation*}
  \widetilde{\mathbf C}^i = \left\{
\begin{array}{rl}
\mathbf{C}^{ev} & \text{ if } i \text{ is even} , \\
\mathbf{C}^{od} & \text{ if }  i \text{ is odd} , 
\end{array} \right.
\end{equation*}
and where the differentials $\widetilde{\partial}: \widetilde{\mathbf C}^i\longrightarrow \widetilde{\mathbf C}^{i+1}$ are restrictions of $\partial$ for all $i\in\mathbb Z$.

\begin{Definition}
A filtered $\mathbb Z_2$-graded complex is a quadruple $(\mathbf C, \partial, \mathbf F^{ev}_{\ast}, \mathbf F^{od}_{\ast})$ where $(\mathbf C=\mathbf C^{ev}\oplus \mathbf C^{od} , \partial)$ is a $\mathbb Z_2$-graded complex, $\ast\in\mathbb N$, $\mathbf F^{ev}_{\ast}$ is a decreasing exhaustive filtration of $\mathbf C^{ev}$, $\mathbf F^{od}_{\ast}$ is a decreasing exhaustive filtration of $\mathbf C^{od}$ and $\partial(\mathbf F^{ev}_i)\subset \mathbf F^{od}_i$,
$\partial(\mathbf F^{od}_i)\subset \mathbf F^{ev}_i$ for all $i\in\mathbb N$. A morphism of filtered $\mathbb Z_2$-graded complexes is a morphism of $\mathbb Z_2$-graded complexes which preserves the filtrations.
\end{Definition}

Given a filtered $\mathbb Z_2$-graded complex $(\mathbf C, \partial, \mathbf F^{ev}_{\ast}, \mathbf F^{od}_{\ast})$,
the associated complex $\widetilde{\mathbf C}$ is filtered by the subcomplexes $\widetilde{\mathbf F}_p$, $p\in\mathbb N$, defined by
\begin{align*}
 \widetilde{\mathbf F}_p^j = \left\{
\begin{array}{rl}
\mathbf F_p^{ev} & \text{ if } j \text{ is even} , \\
\mathbf F_p^{od} & \text{ if }  j \text{ is odd} .
\end{array} \right.
\end{align*}
\vspace{.1mm}

\begin{Proposition}
Let $(\mathbf C, \partial, \mathbf F^{ev}_{\ast}, \mathbf F^{od}_{\ast})$ be a filtered $\mathbb Z_2$-graded complex. Then, there is a spectral sequence $\{E^{p,q}_r, d_r\}_r$ 
with
$
 E_1^{p,q}= H^{p+q}( \widetilde{\mathbf F}_p/  \widetilde{\mathbf F}_{p+1})
$.
If the filtrations $\mathbf F^{ev}_{\ast}, \mathbf F^{od}_{\ast}$ are bounded, this spectral sequence converges to $H^{\bullet}(\widetilde {\mathbf C}, \widetilde{\partial})$, that is, for all $k\in\mathbb Z$ there is a filtration $F_{k,\ast}$ of $H^{k}(\widetilde {\mathbf C}, \widetilde{\partial})$ so that for all $p\in \mathbb N, q\in\mathbb Z$,
\[
 E^{p,q}_{\infty}\cong F_{p+q,\,p}\, H^{p+q}(\widetilde {\mathbf C}, \widetilde{\partial})/ F_{p+q,\,p+1}\, H^{p+q}(\widetilde {\mathbf C}, \widetilde{\partial}).
\]
\end{Proposition}
\begin{Proof}
The spectral sequence is the one associated to the filtered differential complex $(\widetilde{\mathbf C}, \widetilde{F})$, the statement follows e.g. from \cite{MCC}*{Theorem 3.2} or \cite{BT}*{14.6}.
 $\Box$
\end{Proof}
\vs

{\it Twisted cohomology:} Let $(\mathfrak{a}=\oplus_{i\geqslant 0}\mathfrak{a}^{i}, \partial_{\mathfrak{a}})$ be a commutative differential graded  $\mathbf k$--algebra\footnote{Commutativity is understood in the graded sense. From now on, we abbreviate ``differential graded'' by ``dg''.}. Let $(\mathfrak{m}=\oplus_{i\geqslant 0}\mathfrak{m}^{i}, \partial)$ be a dg--$\mathfrak{a}$--module, assume $\mathfrak{m}^{i}=0$ for all but a finite number of
$i\in \mathbb N$. Let 
$t=\sum_{i\geqslant 1}t_{2i+1}\in\mathfrak{a}$ be an element with all components of odd degree greater or equal than three and such that $\partial_{\mathfrak a} t=0$,
set $\partial_t=\partial+ t\,\cdot: \mathfrak{m}\longrightarrow \mathfrak{m}$. If we put $\mathfrak{m}^{od}=\oplus_{i\, odd\,}\mathfrak{m}^i$, $\mathfrak{m}^{ev}=
\oplus_{i\, even\,}\mathfrak{m}^i$, then  $(\mathfrak m=\mathfrak m^{ev}\oplus\mathfrak m^{od}, \partial_t)$ is a $\mathbb Z_2$-graded complex and so we can consider its twisted cohomology. The following proposition is stated in several articles, and it is proved  in \cite{LLW}. The proof we include here avoids the use of exact couples. \vspace{1mm}

\begin{Proposition}
Under the assumptions and with the notations above:
\begin{enumerate}
\item There is a spectral sequence $\{E^{p,q}_r, d_r\}_r$ 
with
\[
 E_2^{p,q}= \left\{
\begin{array}{ll}
H^p(\mathfrak{m},\partial)\ \ \  & \text{if $q$ is even} , \\
\ \  0\ & \text{if $q$ is odd} ,
\end{array} \right.
\]
which converges to the twisted cohomology of $(\mathfrak m, \partial_t)$, in the sense that there are finite decreasing exhaustive filtrations $\{G_i^{ev}\}_{i\geqslant 0}$
of $H^{ev}(\mathfrak m,\partial_t)$ and $\{G_i^{od}\}_{i\geqslant 0}$
of $H^{od}(\mathfrak m, \partial_t)$ such that
\begin{eqnarray*}
E^{p,q}_{\infty} &=& G^{ev}_p H^{ev}(\mathfrak m, \partial_t)/G^{ev}_{p+1} H^{ev}(\mathfrak m, \partial_t) \ \ \ \text{ for }\ p+q \  \text{ even},\\ 
E^{p,q}_{\infty} &=& G^{od}_p H^{od}(\mathfrak m, \partial_t)/G^{od}_{p+1} H^{od}(\mathfrak m, \partial_t) \ \ \ \text{ for }\
p+q \  \text{ odd}.
\end{eqnarray*}
In the sequel, if we want to make explicit the data which gives rise to this spectral sequence, we will denote it by $\{E(\mathfrak m, \partial, t)_r, d_r\}_r$.
\item If $H^{\bullet}(\mathfrak m, \partial)$ is finite dimensional, there is a canonical, functorial isomorphism
\[
\det H^{\bullet}(\mathfrak m, \partial) = \det H^{\bullet}(\mathfrak m, \partial _t) .
\]
\end{enumerate}

\end{Proposition}
\begin{Proof} (see also \cite{LLW}*{3.2}): Setting 
\begin{eqnarray*} 
\mathbf F^{ev}_i&=&\oplus_{\substack{j\geqslant i \\  j \ \text{even}}}
 \mathfrak{m}^{j} \\
\mathbf F^{od}_i&=&\oplus_{\substack{j\geqslant i \\  j \ \text{odd}}}
 \mathfrak{m}^{j}\,,
\end{eqnarray*}
we have a filtered $\mathbb Z_2$-graded complex $(\mathfrak m, \partial_t, \mathbf F^{ev}_{\ast}, \mathbf F^{od}_{\ast})$.
By proposition 1 there is a spectral sequence whose  $E_1$--term is
 \begin{align*}
E_1^{p,q}=H^{p+q}( \widetilde{\mathbf F}_p/  \widetilde{\mathbf F}_{p+1}) = \left\{
\begin{array}{ll}
\mathfrak m^p & \text{ if }  q \text{ is even} , \\
0 & \text{ if }  q \text{ is odd} .
\end{array} \right.
\end{align*}
One checks that $d_1=\partial$, and the filtrations $\mathbf F^{ev}_{\ast}, \mathbf F^{od}_{\ast}$ are clearly bounded, so  i) follows from  proposition 1. In particular, if $H^{\bullet}(\mathfrak m, \partial)$ is finite dimensional, then so is $H^{\bullet}(\mathfrak m, \partial_t)$. \vs

For ii), notice that for each $r\geqslant 2$ and $p,q\in\mathbb Z$, we have $E^{p,q}_r=E^{p,q+2}_r$, so if we put 
\[
 E^{ev}_r=\oplus_{p+q=0} E^{p,q}_r \ \ , \ \ E^{od}_r=\oplus_{p+q=1} E^{p,q}_r , 
\]
then we have a $\mathbb Z_2$-graded complex
$( E^{ev}_r\oplus E^{od}_r, d_r)$ and $H^{\bullet}( E^{ev}_r\oplus E^{od}_r, d_r)=E^{ev}_{r+1}\oplus E^{od}_{r+1}$. Then, one obtains the isomorphism in ii) from i) and lemma 1, functoriality is easily checked. $\Box$
\end{Proof}\vspace{0.1cm}

The following elementary lemma will be useful in the sequel:

\begin{Lemma}
For $i=1,2$, let $\mathfrak a_i=\oplus_{j\geqslant 0}\mathfrak a^j_i$ be commutative dg-$\mathbf k$-algebras, let $(\mathfrak m_i=\oplus_{j\geqslant 0}\mathfrak m^j_i, \partial_i)$ be dg-modules over 
$\mathfrak a_i$, assume $\mathfrak m^j_1=\mathfrak m^j_2=0$ for all but a finite number of indexes $j\geqslant 0$. Let $t_i\in\mathfrak a_i$ be closed elements with all components of odd degree greater or equal than three, set $\partial_{t_i}=\partial_i + t_i\cdot$. Let $f: (\mathfrak m_1,\partial_1) \longrightarrow (\mathfrak m_2,\partial_2)$  be a  morphism of  $\mathbb Z$-graded complexes which is also a morphism of $\mathbb Z_2$-graded complexes $f_{tw}: (\mathfrak m_1,\partial_{t_1}) \longrightarrow (\mathfrak m_2,\partial_{t_2})$\footnote{We denote $f$ by $f_{tw}$ when we regard it as a morphism between the twisted complexes.}.
%
If $f$ is a quasi-isomorphism, then $f_{tw}$ induces isomorphisms on twisted cohomology.
\end{Lemma}

\begin{Proof}\
Under the assumptions of the lemma, we have a morphism of filtered $\mathbb Z_2$-graded complexes
$f_{tw}: (\mathfrak m_1, \partial_{t_1}, \mathbf F^{ev}_{\ast}, \mathbf F^{od}_{\ast})\longrightarrow (\mathfrak m_2, \partial_{t_2}, \mathbf F^{ev}_{\ast}, \mathbf F^{od}_{\ast})$, and so  an induced morphism $\{f_n\}_{n\geqslant 1}: E(\mathfrak m_1, \partial_1, t_1)_n \longrightarrow E(\mathfrak m_2, \partial_2, t_2)_n$ between the spectral sequences introduced in the previous proposition. The morphism between their abutments is the one induced by $f_{tw}$ on twisted cohomology (see e.g. \cite{MCC}*{Theorem 3.5}). The rows of $E(\mathfrak m_i, \partial_i, t_i)^{p,q}_1$ are zero for $q$ odd and are equal to the untwisted complex $(\mathfrak m_i,\partial_i)$ for $q$ even ($i=1,2$), the map $f_1$ is given by $f$ between the even numbered rows of the $E_1$-terms. If $f$ is a quasi-isomorphism, the $E_2$-terms are isomorphic via $f_2$, and then by standard arguments so are the abutments. $\Box$
\end{Proof}\vs

{\it Reidemeister torsion} (see \cite{Mu}*{section 1}, \cite{MW}*{5.2}):  Given a simplicial complex $K$, we denote by $C_{\bullet}(K)$ the  complex of simplicial chains of $K$ with real coefficients. Assume that $K$ is finite and connected, let $\pi: \widetilde{K}\longrightarrow K$ denote the universal covering space of $K$ with the induced simplicial structure, orient the simplexes of $K$ and $\widetilde{K}$ so that the restriction of $\pi$ to each simplex is orientation preserving. Notice that the action of $\pi_1(K)$ by covering transformations endows $C_{\bullet}(\widetilde{K})$ with a structure of right $\mathbb R[\pi_1(K)]$--module.\vs

Let $\rho: \pi_1(K) \longrightarrow GL(E)$ denote a real or complex finite dimensional representation of the fundamental group of $K$, denote by $H^{\bullet}(K,E)$ the cohomology of the complex
$
 C^{\bullet}(K,E):= Hom_{R[\pi_1(K)]}(C_{\bullet}(\widetilde{K}), E).
$
Choose an embedding of $K$ as a fundamental domain in $\widetilde{K}$, an order on the set of $i$-simplexes of $K$,
a metric on $E$ (hermitian in the complex case) and a unit volume element attached to the metric. These choices determine volume elements $\mu_i\in \det C^i(K,E)$, and the 
\emph{Reidemeister torsion}
$\tau(K, E)\in \det H^{\bullet}(K,E)$ is defined as the image of $\otimes_i\, \mu_i^{(-1)^i}\in \det C^{\bullet}(K,E)$ under the Knudsen-Mumford isomorphism 
$\det C^{\bullet}(K,E) \cong \det H^{\bullet}(K,E)$. 
If the representation 
$\rho$ is unimodular (that is, if the norm of $\det(\rho(g))$ is one for all $g\in \pi_1(K)$), then $\tau(K, E)$
is independent of the choice of the embedding of $K$ into $\widetilde{K}$, up to a scalar of unit norm. If $K$ is homeomorphic to a compact, oriented, odd dimensional manifold, then the Reidemeister torsion is independent of the chosen metric on $E$ (again, up to a scalar of unit norm)\footnote{One can define Reidemeister torsion using densities instead of volume forms, then this indeterminacy is built--in the definition. Alternatively, one can define Reidemeister torsion as a norm on the determinant line.}. \vs

{\it Dupont's simplicial currents} (see \cite{Du}): 
For $q\in\mathbb N$, denote $\Delta^q$ the standard $q$-dimensional simplex and 
$\varepsilon^i :\Delta^{q-1} \longrightarrow \Delta^q$ the inclusions of its faces ($0\leqslant i \leqslant q$).
Let $K$ be a simplicial complex, $K_q$ the set of its $q$--dimensional simplexes, $\varepsilon_i: K_q \longrightarrow K_{q-1}$ the face operators. \vs

Denote by
$A^{\bullet}(K)=\oplus_{k\in\mathbb N} A^k(K)$ 
the Sullivan--Thom--Whitney dg--algebra of piecewise smooth real differential forms on $K$.
Recall that elements of $A^{k}(\Delta^q)$ are smooth $k$--differential forms on $\Delta^q$, and elements of $A^{k}(K)$ are families $\{\varphi_{\sigma}\}_{\sigma}$, where $\sigma\in \bigcup_{q} K_q$, such that
\begin{enumerate}
\item $\varphi_{\sigma}\in A^k (\Delta^q)$ \  if\ $\sigma\in K_q$, 
\item  $\varphi_{\varepsilon_i(\sigma)}=(\varepsilon^i)^{\ast}(\varphi_{\sigma})$ \ for all $q\geqslant 0$, $0\leqslant i\leqslant q$, $\sigma\in K_q$.
\end{enumerate}
The differential is given by $(d\varphi)_{\sigma}=d\varphi_{\sigma}$,
we refer to \cite{FHT}*{Chapter 10} or \cite{BT}*{Chapter 9} for more information on these algebras and their properties. They were used in the context of twisted cohomology in \cite{AS}*{\S 6}, for the purpose of dealing with topological spaces which are not manifolds. \vs

In \cite{Du}, Dupont introduces a homological complex $Du_{\bullet}(K)$ of piecewise smooth currents on $K$, defined as follows: 
If $n\in\mathbb N$, then
\[
Du_{n}(K)=\bigoplus_{j=0}^{\infty} A^j(\Delta^{n+j})\otimes C_{n+j}(K),
\]
the differential $\partial_n: Du_n \longrightarrow Du_{n-1}$ is defined by
\[
 \partial_n(\omega \otimes \sigma)=(-1)^n \,d\omega \otimes \sigma + \sum_{i=0}^{n+j}(-1)^i (\varepsilon^i)^{\ast}(\omega)\otimes \varepsilon_i (\sigma), 
\]
where $\omega\in A^j(\Delta^{n+j})$ and $\sigma\in K_{n+j}$.
The following properties hold, see \cite{Du}:\vs

\begin{Proposition}\hspace{-0.5cm}
\begin{enumerate}
\item The Dupont complex is functorial with respect to simplicial maps. More precisely,
a simplicial map $f:K \longrightarrow K'$ induces a chain map $f_{\bullet}: Du_{\bullet}(K) \longrightarrow Du_{\bullet}(K')$ which, for $\omega\in A^{k}(\Delta^{n+k})$ and $\sigma\in K_{n+k}$, is given by $\omega\otimes \sigma \longmapsto \omega \otimes f(\sigma)$.
\item The product $A^k(K) \times Du_n(K)   \longrightarrow Du_{n-k}(K)$ defined by 
$$\varphi\cdot (\omega \otimes \sigma)= (\varphi_{\sigma}\wedge \omega)\otimes \sigma$$
endows $Du_{\bullet}(K)$ with a dg-module structure over the dg-algebra $A^{\bullet}(K)$.
\item The chain map $\psi: C_{\bullet}(K) \longrightarrow Du_{\bullet}(K)$ given by $\psi(\sigma)=1\otimes \sigma$ induces an isomorphism in homology. 
\end{enumerate}
\end{Proposition}\vs

By i), the group $\pi_1(K)$ acts on the right on $Du_{\bullet}(\widetilde{K})$ by covering transformations, so we can consider the complex of $\mathbb R$-vector spaces $\mathfrak m^{\bullet}(K,E):=Hom_{\mathbb R[\pi_1(K)]}
(Du_{\bullet}(\widetilde{K}), E)$. It follows from ii) above that $\mathfrak m^{\bullet}(K,E)$ is endowed with a structure of $A^{\bullet}(K)$--dg--module.
\vs

\begin{Definition}
Let $K$ be a finite simplicial complex, let $\rho: \pi_1(K) \longrightarrow Aut(E)$ be a real or complex finite dimensional representation of the fundamental group of $K$, denote by $\partial$ the differential of $\mathfrak m^{\bullet}(K,E)$. Given $t\in \oplus_{i\geqslant 1} A^{2i+1}(K)$ closed, put
\[
\partial_t=\partial + t \cdot: \mathfrak m^{\bullet}(K,E)\longrightarrow \mathfrak m^{\bullet}(K,E)\,.
\]
Then, $(\mathfrak m^{\bullet}(K,E), \partial_t)$ is a $\mathbb Z_2$--graded complex, and we define the \emph{twisted simplicial cohomology} attached to $K,E,t$ by
\[
 H^{ev}(K, E, t):=H^{ev}(\mathfrak m^{\bullet}(K,E), \partial_t)\ \ \ \text{,}\ \ \ H^{od}(K, E, t):=H^{od}(\mathfrak m^{\bullet}(K,E), \partial_t).
\]
We put $\det H^{\bullet}(K, E, t):= \det H^{\bullet}(\mathfrak m^{\bullet}(K,E), \partial_t)$. 
\end{Definition}

\begin{Remarks}
i) By definition, $H^{ev}(K, E, t)$ and $H^{od}(K, E, t)$ are real vector spaces. However, if $\rho$ is a rational representation and, instead of piecewise smooth forms, one considers piecewise polynomial forms with rational coefficients, then,
in the same way as above, one can define a twisted simplicial cohomology over $\mathbb Q$. \vs
 
ii) Twisted simplicial cohomology (and similarly twisted de Rham cohomology) depends on the piecewise smooth form $t$ and not just on the cohomology class it defines. Still, if we assume for simplicity that $t$ is of pure odd degree $2i+1$, we can put
\[
\widetilde{H}^{\bullet}(K,E, [t])=\varinjlim_{b\in A^{2i}(K)}H^{\bullet}(K,E, t+ \partial b)\,,
\]
where $\bullet=ev, od$ and the direct system is given by the transition maps 
\begin{eqnarray*}
H^{\bullet}(K,E, t+ \partial b_1) &\longrightarrow& H^{\bullet} (K,E, t+ \partial b_2).  \\ 
{[}\omega{]} &\longmapsto& {[}e^{b_1-b_2}\omega{]}
\end{eqnarray*}
In fact, since all transition maps are isomorphisms, one could work with inverse limits as well. By its very definition, 
$\widetilde{H}^{\bullet}(K,E, [t])$ depends only on the cohomology class defined by $t$.\vs

We prove next functoriality and invariance under subdivision of $H^{\bullet}(K, E, t)$:
\end{Remarks}
\vspace{2mm}


\begin{Proposition}
\begin{enumerate}
\item Let $K, L$ be simplicial complexes, $\rho: \pi_1(L) \longrightarrow Aut(E)$ a finite dimensional linear representation,
$t_L\in \oplus_{i\geqslant 1} A^{2i+1}(K)$ such that $\partial\,t_L=0$, $f: K \longrightarrow L$ a simplicial map, $f: A^{\bullet}(L)\longrightarrow A^{\bullet}(K)$ the induced morphism of dg-algebras, $t_K=f(t_L)$. 
Then, $f$ induces linear maps 
\[
 f^{\bullet}:  H^{\bullet}(L, E, t_L) \longrightarrow H^{\bullet}(K, f^{\ast}E, t_K)
\]
for $\bullet=ev, od$, where $f^{\ast}E$ denotes the representation \\ $\pi_1(K) \stackrel{f_{\ast}}\longrightarrow \pi_1(L) \longrightarrow Aut(E)$.
\item If $g: J \longrightarrow K$ is another simplicial map, then  $(f\circ g)^{\bullet}=g^{\bullet}\circ f^{\bullet}$,
for $\bullet=ev, od$, that is, twisted cohomology is a contravariant functor.
\end{enumerate}
\end{Proposition}
\begin{Proof}
It follows from the definitions that the given data  induce a morphism of filtered $\mathbb Z_2$-graded complexes
$
(\mathfrak m^{\bullet}(L,E), \partial_{t_L}) \longrightarrow  
(\mathfrak m^{\bullet}(K,f^{\ast}E),\partial_{t_K})
$, the maps $f^{\bullet}$ are the induced maps in cohomology, the proof of ii) is straightforward.
$\Box$
\end{Proof}

\begin{Proposition}
 Let $K$ be a simplicial complex, $\rho: \pi_1(K) \longrightarrow Aut(E)$ a finite dimensional linear representation of its fundamental group, 
$t\in \oplus_{i\geqslant 1} A^{2i+1}(K)$ closed. If $K'$ is a subdivision of $K$, restriction induces a natural quasi-isomorphism of dg-algebras $\varsigma: A^{\bullet}(K)\longrightarrow A^{\bullet}(K')$. If $t'=\varsigma(t)$, then 
there are canonical isomorphisms
\[
 H^{ev}(K, E, t) \cong H^{ev}(K', E, t') \ \ \text{ and }\ \ H^{od}(K, E, t)\cong  H^{od}(K', E, t').
\]
\end{Proposition}

\begin{Proof} By functoriality of twisted cohomology and the simplicial approximation theorem 
(see e.g. \cite{Mun}*{Chapter 2}), we can assume there is a map $g: K' \longrightarrow K$ which is a simplicial 
approximation of the identity map between the underlying topological spaces. Then $g$ induces a map of complexes $\widetilde{g}: (\mathfrak m^{\bullet}(K,E), \partial) \longrightarrow  
(\mathfrak m^{\bullet}(K',E), \partial)$ which is a quasi--isomorphism  by invariance of cohomology under subdivision. The map $g$ induces also a 
morphism of  
$\mathbb Z_2$-graded complexes 
$
\widetilde{g}_{tw}:(\mathfrak m^{\bullet}(K,E), \partial_{t}) \longrightarrow  
(\mathfrak m^{\bullet}(K',E),\partial_{t'})
$.
By lemma 2, the map induced by $\widetilde{g}_{tw}$ on twisted cohomology is an isomorphism, and so we are done. $\Box$
\end{Proof} \vspace{1.5mm}

\begin{Lemma} 
i) The chain map 
$
\psi^{\ast}: (\mathfrak m^{\bullet}(K,E), \partial)\longrightarrow (C^{\bullet}(K,E), \partial)
$
induced by $\psi: C_{\bullet}(\widetilde{K}) \longrightarrow Du_{\bullet}(\widetilde{K})$ is a quasi-isomorphism. \vspace{1mm}

ii) There is a canonical isomorphism
\[
 \det H^{\bullet}(K, E)
= \det H^{\bullet}(K, E, t).
\]
\end{Lemma}
\begin{Proof}
Both $C_{\bullet}(\widetilde{K})$ and  $Du_{\bullet}(\widetilde{K})$  are complexes of free $\mathbb R[\pi_1(K)]$-modules, a basis of
$C_{n}(\widetilde{K})$  is given by the elements of $ K_n$,  regarded as elements of $C_n(\widetilde{K})$ via the choice of an inclusion of $K$ as a fundamental domain in $\widetilde{K}$, a basis of $ Du_{n}(\widetilde{K})$ is given by the elements of the form $\varphi\otimes \sigma$, where $\sigma\in K_{n+j}$, and $\varphi\in A^{j}(\Delta^{n+j})$, $j\geqslant 0$ runs over a basis of $A^{j}(\Delta^{n+j})$ as an $\mathbb R$--vector space. Since $\psi$ is a quasi-isomorphism, it follows that $\psi^{\ast}$ is a quasi-isomorphism as well, as claimed in i). \vs

To prove ii), notice that from i) follows that 
\[
\det H^{\bullet}(K, E):=\det (C^{\bullet}(K,E), \partial)=\det H^{\bullet}(\mathfrak m^{\bullet}(K,E), \partial).
\] 
By proposition 2, we also have 
\[
\det H^{\bullet}(\mathfrak m^{\bullet}(K,E), \partial) = \det H^{\bullet}(\mathfrak m^{\bullet}(K,E), \partial_t)=\det H^{\bullet}(K, E, t),
\] 
and then the claimed isomorphism follows.
$\Box$
\end{Proof}

\begin{Definition} 
Assume that the representation $\rho$ is unimodular.  We define the \emph{twisted Reidemeister torsion} attached to $K,E,t$, denoted
\[
\tau_{twist} (K, E, t)\in \det H^{\bullet}(K, E, t)
\] 
as the image of the Reidemeister torsion $\tau(K, E)\in \det H^{\bullet}(K,E)$, under the isomorphism in lemma 3.ii). The following invariance under subdivision holds:
\end{Definition}

\begin{Proposition} Let $K$ be a simplicial complex, $\rho: \pi_1(K) \longrightarrow Aut(E)$ a unimodular representation of its fundamental group, $t\in \oplus_{i\geqslant 1} A^{2i+1}(K)$ with $\partial\,t=0$. Let $K'$ be a subdivision of $K$, $\varsigma: A^{\bullet}(K)\longrightarrow A^{\bullet}(K')$ the restriction map, $t'=\varsigma(t)$. Denote
\[
\iota: \det H^{\bullet}(K, E, t) =  \det H^{\bullet}(K', E, t')
\]
the isomorphism induced by the isomorphisms in proposition 5. Then, up to a scalar of norm one, 
$\iota(\tau_{twist} (K, E, t)) = \tau_{twist} (K', E, t')$. 
 \end{Proposition}
\begin{Proof} By lemma 3.i), via the isomorphism $\psi^{\ast}$ we can view  Reidemeister torsion as an element of\, $\det H^{\bullet}(\mathfrak m^{\bullet}(K,E), \partial)$. 
As in the proof of proposition 5, we can assume we have a simplicial approximation of the identity map $K' \longrightarrow K$. By functoriality in proposition 2.ii), we have a commutative square of isomorphisms
 \[
\begin{CD}
{\det H^{\bullet}(\mathfrak m^{\bullet}(K,E), \partial)} @>f_K>> {\det H^{\bullet}(\mathfrak m^{\bullet}(K,E), \partial_t)} \\
@ VVV @ VV\iota V \\
{\det H^{\bullet}(\mathfrak m^{\bullet}(K',E), \partial)} @>>f_{K'}> {\det H^{\bullet}(\mathfrak m^{\bullet}(K',E), \partial_{t'})} \\
\end{CD}
\] 
where the horizontal arrows are  given by proposition 2. 
It is known that, up to a scalar of unit norm, (untwisted) Reidemeister torsion is invariant under subdivision.
The twisted Reidemeister torsion attached to $K,E,t$ (respectively, $K',E,t'$) is the image of the untwisted one under $f_K$ (resp., $f_{K'}$). Then, the result follows. $\Box$\vs
\end{Proof}

\begin{Remark}
In \cite{AAS}, C. Arias Abad and F. Sch\"atz defined a Reidemeister torsion for flat superconnections. Their construction can be applied to $(\mathfrak m^{\bullet}(K,E),\partial_{t})$ and, as F. Sch\"atz explained to me, their torsion coincides (up to elements of unit norm) with the one defined above. We also remark that V. Mathai and S. Wu defined in \cite{MW2} a Reidemeister torsion for $\mathbb Z_2$--graded elliptic complexes.
\end{Remark}
\vs

{\bf Comparison of twisted cohomologies and torsions.}\vs

{\it Comparison of twisted cohomologies in the smooth case:} Let $X$ be a compact, connected, oriented smooth manifold,
denote by
$\Omega^{\bullet}(X)$ the de Rham algebra of smooth differential forms on $X$, denote by $H^{\bullet}(X, \mathbb R)$ its cohomology, let $T=\sum_{i\geqslant 1}T_{2i+1}\in\Omega^{\bullet}(X)$ be a closed differential form on $X$ with all components of odd degree greater or equal than three.  Let $(\mathcal E, \nabla)$ be a flat vector bundle over $X$, put
$\Omega^{\bullet}(X, \mathcal E)=(\Omega^{\bullet}\otimes \mathcal E)(X)$. The connection $\nabla: \mathcal E(X) \longrightarrow (\Omega^1\otimes \mathcal E)(X)$ extends to a  differential $\delta: \Omega^{\bullet}(X, \mathcal E)\longrightarrow \Omega^{\bullet+1}(X, \mathcal E)$,
denote $\delta_T$ the twisted differential $\delta+ T\,\wedge\cdot$\,. 
Following \cite{MW}, \cite{RW},  define the twisted de Rham cohomology 
groups attached to the data $X, \mathcal E, T$ as $H^{ev}(X,\mathcal E, T):=H^{ev}(\Omega^{\bullet}(X, \mathcal E), \delta_T)$, $H^{od}(X,\mathcal E, T):=H^{od}(\Omega^{\bullet}(X, \mathcal E), \delta_T)$.\vs

Let $\mathcal E^{\nabla}$ denote the local system of flat sections of $\mathcal E$, $\pi:\widetilde{X}\longrightarrow X$ the universal cover of $X$,  $E=\Gamma(\widetilde{X},\mathcal \pi^{\ast}\mathcal E^{\nabla})$ 
the space of global sections of the pull-back of $\mathcal E^{\nabla}$, denote $\rho: \pi_1(X) \longrightarrow Aut(E)$ the monodromy representation. 
Endow $X$ with a smooth triangulation given by a finite simplicial complex $K$ and $\widetilde{X}$ with the induced triangulation. Then, we have the following comparison of twisted cohomologies (see also \cite{MW}*{Lemma 5.3}): 

\begin{Proposition}
There are canonical isomorphisms
\[
 H^{ev}(K, E, T)\cong H^{ev}(X,\mathcal E, T) \ \ , \ \ 
 H^{od}(K, E, T)\cong H^{od}(X,\mathcal E, T)\,.
\]
\end{Proposition}

\begin{Proof} 
Given  $i,n\geqslant 0$, let $\widetilde{\sigma}$ be a $(n+i)$-dimensional simplex in $\widetilde{K}$, let $\sigma=\pi(\widetilde{\sigma})$. If $\omega\in \Omega^n(X,\mathcal E)$, in a small neighborhood $U$ of  $\sigma$ we can write $\omega=\sum \omega_i \otimes s_i$ where 
$\omega_i\in \Omega^n(X)$ and $s_i$ is a local section of $\mathcal E^{\nabla}$. Let $\widetilde{U}$ be the connected component of $\pi^{-1}(U)$ containing $\widetilde{\sigma}$. Then, $s_i\circ \pi_{\mid \widetilde{U}}$ can be identified with an element of $E$, for simplicity we denote it also by $s_i$. For $n\geqslant 0$, define
\begin{eqnarray*}
\Omega^n(X,\mathcal E) &\stackrel{\mathfrak I^n}\longrightarrow& \mathfrak m^{n}(K,E)=Hom_{\mathbb R[\pi_1(K)]} \left(\bigoplus_{i=0}^{\infty} A^i(\Delta^{n+i})\otimes C_{n+i}(\widetilde{K}), E \right). \\
\omega \ \ \ &\longmapsto& \ \ \ \ \ \ \ \left[ \varphi \otimes \sigma \mapsto \sum_i \left( \int_{\sigma}  \omega_i \wedge \varphi \right) s_i \right]
\end{eqnarray*}
It is easily checked that this defines both a morphism of $\mathbb Z$-graded complexes  $\mathfrak I^{\bullet}: (\Omega^{\bullet}(X,\mathcal E), \delta)\longrightarrow (\mathfrak m^{\bullet}(K,E), \partial)$ and a morphism  of $\mathbb Z_2$-graded complexes 
$\mathfrak I^{\bullet}_{tw}: 
(\Omega^{\bullet}(X,\mathcal E), \delta_T) \longrightarrow  (\mathfrak m^{\bullet}(K,E), \partial_T)$. 
Let $dR: (\Omega^{\bullet}(X,\mathcal E), \delta) \longrightarrow (C^{\bullet}(K,E), \partial)$ be the integration map which gives the de Rham theorem with coefficients (as described for example in \cite{Mu}*{section 1}). 
We have $dR= \psi^{\ast}\circ \mathfrak I^{\bullet}$, where $\psi^{\ast}$ is induced by the map 
$\psi: C_{\bullet}(\widetilde{K})\longrightarrow Du_{\bullet}(\widetilde{K})$ in proposition 3, iii). Since both $dR$ and $\psi^{\ast}$ are quasi-isomorphisms (for the first, see \cite{Mu}), $\mathfrak I^{\bullet}$ is a quasi-isomorphism as well and then, by lemma 2, $\mathfrak I^{\bullet}_{tw}$ induces an isomorphism on cohomology.
$\Box$
\end{Proof}\vs

{\it Comparison with Mathai--Wu's definition:}  
Let $K$ be a finite connected simplicial complex, let  $\rho: \pi_1(K) \longrightarrow GL(E)$ be a unimodular representation of its fundamental group, let
$\vartheta=\sum_{i\geqslant 0}\vartheta_{2i+1}\in C^{\bullet}(K)$ be a closed cochain with all its components of odd degree greater than $\dim K/2$ (here, $\dim K$ is the maximum of the dimensions of its simplexes).
%
Then $\vartheta \cup\vartheta =0$, so one has a $\mathbb Z_2$--graded cochain complex  
$(C^{\bullet}(K,E), \partial_\vartheta)$, where $\partial_\vartheta= \partial +\vartheta\,\cup\cdot $. Set
\[
 H_{MW}^{ev}(K, E, \vartheta):= H^{ev}(C^{\bullet}(K,E), \partial_{\vartheta}) 
\ , \ H_{MW}^{od}(K, E, \vartheta):= H^{od}(C^{\bullet}(K,E), \partial_{\vartheta}).
\]
In \cite{MW}*{5.2}, the authors define  
a twisted Reidemeister torsion $\tau_{MW}(K,E,\vartheta)\in \det H_{MW}^{\bullet}(K,E,\vartheta)$ as the image of $\otimes_i\, \mu_i^{(-1)^i}\in \det C^{\bullet}(K,E)$ under the isomorphism $\det C^{\bullet}(K,E) \cong \det 
H^{\bullet}_{MW}(K,E,\vartheta)$. We have:\vs

\begin{Proposition}
With the notations above, for each $i\geqslant 0$ let $t_{2i+1}\in A^{2i+1}(K)$ be a closed preimage of $\vartheta_{2i+1}\in C^{2i+1}(K)$
under the integration map \footnote{Such preimages always exist, because the integration map is a quasi-isomorphism.}
\begin{eqnarray*}
 A^{2i+1}(K) &\longrightarrow & C^{2i+1}(K).\\
\omega &\longmapsto & \left[\sigma \mapsto \int_{\sigma} \omega_{\sigma}\right]
\end{eqnarray*}
Set $t=\sum_{i \geqslant 0}t_{2i+1}$. Then, there are canonical isomorphisms
\[
  H_{MW}^{ev}(K,E,\vartheta) \cong  H^{ev}(K,E, t) \ \ \text{and}\ \   H_{MW}^{od}(K,E,\vartheta) \cong  H^{od}(K,E, t) ,
\]
which induce an isomorphism
\[
\iota: \det H_{MW}^{\bullet}(K,E,\vartheta) \longrightarrow \det H^{\bullet}(K,E, t)
\]

such that $\iota (\tau_{MW} (K, E, \vartheta)) = \tau_{twist}(K, E, t)$. 
\end{Proposition}
\begin{Proof} 
 The map $\psi$ in proposition 3.iii) induces a morphism of filtered $\mathbb Z$-graded complexes  $\psi^{\ast}:(\mathfrak m^{\bullet}(K,E),\partial) \longrightarrow  (C^{\bullet}(K,E), \partial)$ and also a morphism of $\mathbb Z_2$--graded complexes $\psi^{\ast}_{tw}:(\mathfrak m^{\bullet}(K,E),\partial_t) \longrightarrow  (C^{\bullet}(K,E), \partial_{\vartheta})$. The former is a quasi-isomorphism by lemma 3.i), so, by lemma 2, the morphism $\psi^{\ast}_{tw}$ induces an isomorphism between twisted cohomologies. \vs

The induced isomorphism between determinant lines factors as the composition 
\begin{eqnarray*}
\det H^{\bullet}(C^{\bullet}(K,E), \partial_{\vartheta}) = \det H^{\bullet}(C^{\bullet}(K,E), \partial) = \\
\det H^{\bullet}(\mathfrak m^{\bullet}(K,E), \partial) = \det H^{\bullet}(\mathfrak m^{\bullet}(K,E), \partial_t). \ \ \ \ \
\end{eqnarray*}
Here, the first and third equalities (that is, canonical isomorphisms) follow from proposition 2, ii), while the second is induced by the map 
$\psi^{\ast}$ considered above.
Then, the statement about twisted Reidemeister torsions follows from the definitions. $\Box$
\end{Proof}

\bibliographystyle{amsplain}
\bibliography{mybibliography.bib}

%
%
%
%
%
\end{document}